\documentclass[english]{article}
\usepackage[T1]{fontenc}
\usepackage[latin9]{inputenc}
\usepackage{geometry}
\usepackage[active]{srcltx}
\usepackage{amsmath}
\usepackage{amsthm}
\usepackage{amssymb}
\usepackage{graphicx}

\makeatletter
\numberwithin{equation}{section}
\numberwithin{figure}{section}

\usepackage{color}
\usepackage{amsthm}
\theoremstyle{plain}

\usepackage{hyperref}

\makeatother

\usepackage{babel}
\begin{document}
\global\long\def\ve{\varepsilon}%
\global\long\def\R{\mathbb{R}}%
\global\long\def\Rn{\mathbb{R}^{n}}%
\global\long\def\Rd{\mathbb{R}^{d}}%
\global\long\def\E{\mathbb{E}}%
\global\long\def\P{\mathbb{P}}%
\global\long\def\bx{\mathbf{x}}%
\global\long\def\vp{\varphi}%
\global\long\def\ra{\rightarrow}%
\global\long\def\smooth{C^{\infty}}%
\global\long\def\Tr{\mathrm{Tr}}%
\global\long\def\bra#1{\left\langle #1\right|}%
\global\long\def\ket#1{\left|#1\right\rangle }%
\global\long\def\Re{\mathrm{Re}}%
\global\long\def\Im{\mathrm{Im}}%
\global\long\def\bsig{\boldsymbol{\sigma}}%
\global\long\def\btau{\boldsymbol{\tau}}%
\global\long\def\bmu{\boldsymbol{\mu}}%
\global\long\def\bx{\boldsymbol{x}}%
\global\long\def\bups{\boldsymbol{\upsilon}}%
\global\long\def\bSig{\boldsymbol{\Sigma}}%
\global\long\def\bt{\boldsymbol{t}}%
\global\long\def\bs{\boldsymbol{s}}%
\global\long\def\by{\boldsymbol{y}}%
\global\long\def\brho{\boldsymbol{\rho}}%
\global\long\def\ba{\boldsymbol{a}}%
\global\long\def\bb{\boldsymbol{b}}%
\global\long\def\bz{\boldsymbol{z}}%
\global\long\def\bc{\boldsymbol{c}}%
\global\long\def\balpha{\boldsymbol{\alpha}}%
\global\long\def\bbeta{\boldsymbol{\beta}}%
\global\long\def\blam{\boldsymbol{\lambda}}%
\global\long\def\Blam{\boldsymbol{\Lambda}}%

\title{MNE: overparametrized neural evolution with \\
applications to diffusion processes and sampling}
\author{Michael Lindsey\vspace{1mm}\\
{\normalsize{}UC Berkeley and LBNL}}
\date{}
\maketitle
\begin{abstract}
We propose a framework for solving evolution equations within parametric
function classes, especially ones that are specified by neural networks.
We call this framework the minimal neural evolution (MNE) because
it is motivated by the goal of seeking the smallest instantaneous
change in the neural network parameters that is compatible with exact
solution of the evolution equation at a set of evolving collocation
points. Formally, the MNE is quite similar to the recently introduced
Neural Galerkin framework, but a difference in perspective motivates
an alternative sketching procedure that effectively reduces the linear
systems solved within the integrator to a size that is interpretable
as an effective rank of the evolving neural tangent kernel, while
maintaining a smooth evolution equation for the neural network parameters.
We focus specifically on the application of this framework to diffusion
processes, where the score function allows us to define intuitive
dynamics for the collocation points. These can in turn be propagated
jointly with the neural network parameters using a high-order adaptive
integrator. In particular, we demonstrate how the Ornstein-Uhlenbeck
diffusion process can be used for the task of sampling from a probability
distribution given a formula for the density but no training data.
This framework extends naturally to allow for conditional sampling
and marginalization, and we show how to systematically remove the
sampling bias due to parametric approximation error. We validate the
efficiency, systematic improvability, and scalability of our approach
on illustrative examples in low and high spatial dimensions.
\end{abstract}

\section{Introduction \label{sec:Introduction}}

Neural network-based methods for solving differential equations, in
both low and high spatial dimensions, are gaining increasing attention
due to the capacity of neural network architectures to flexibly represent
complicated functions. In particular, several frameworks have emerged
for propagating neural network approximations in time. Though we do
not attempt a comprehensive review, we point out that a contrast has
emerged between global-in-time approaches and sequential-in-time approaches.
The former category includes PINNS \cite{Raissi_Perdikaris_Karniadakis_2019}
among many other approaches. But our focus in this work is on the
latter category, which includes Neural Galerkin schemes \cite{Bruna_Peherstorfer_Vanden-Eijnden_2024},
regularized dynamical parametric approximation \cite{feischl2024regularizeddynamicalparametricapproximation},
TENG \cite{TENG}, and evolutional deep neural networks \cite{PhysRevE.104.045303},
to name a few. For more comprehensive reviews, see, e.g., \cite{Bruna_Peherstorfer_Vanden-Eijnden_2024,WEN2024134129,BermanPeherstoferSparse,zhang2024sequentialintimetrainingnonlinearparametrizations}.

Neural Galerkin schemes relate closely to optimization approaches
based on natural gradients \cite{amari1985differential,amari1998natural,doi:10.1137/22M1477805},
including methods for optimizing quantum many-body wavefunctions within
the Variational Monte Carlo (VMC) framework \cite{gubernatis2016quantum,becca2017quantum}.
Wavefunction optimization based on natural gradients has a long history
in the quantum chemistry literature, where it is known traditionally
as stochastic reconfiguration (SR) \cite{becca2017quantum}. Recently,
a simple modification of stochastic reconfiguration called MinSR \cite{ChenHeyl}
and various extensions \cite{GOLDSHLAGER2024113351} have emerged
as leading optimizers in VMC. Our work is motivated in part by the
contrast in perspective between SR and MinSR.

Neural Galerkin and natural gradient approaches can be derived by
projecting the dynamics that are prescribed in the function space
onto the tangent space of the neural parametric manifold. At first
glance, the natural gradient perspective suggests that it is desirable
to take a large number $N$ of sample points---in particular, a number
greater than the number of neural network parameters $p$---in order
to compute this projection accurately. Meanwhile, MinSR can be viewed
as attempting to solve the prescribed dynamics exactly at a set of
points, which in this work we interpret as collocation points. The
philosophy of MinSR is that $p>N$, so that this prescription is underdetermined,
and the \emph{smallest} possible change to the neural network parameters
is chosen among those that solve the dynamics exactly at the collocation
points.

In fact, the SR and MinSR updates are mathematically equivalent, if
the same regularization parameter is used for both methods \cite{GOLDSHLAGER2024113351},
regardless of whether $p>N$ or $p<N$. (For SR, this regularization
constitutes a shift of the metric. For MinSR, it softens the constraint
at the collocation points.) Relatedly, it has been observed that in
Neural Galerkin schemes, the parameter update minimizes a least squares
objective \cite{zhang2024sequentialintimetrainingnonlinearparametrizations},
which corresponds to a soft version of the MinSR constraint. We point
out that the impact of regularization on dynamical parametric approximation,
including the Neural Galerkin framework, has been studied recently
\cite{feischl2024regularizeddynamicalparametricapproximation}, though
the focus of this work does not address the choice of collocation
points.

Within Neural Galerkin schemes, key computational questions arise
about how to deal with the least squares problems called as subroutines
by the integrator, as well as how to integrate the parameter trajectories
while simultaneously propagating the particles (which we will view
as collocation points) to resolve the evolving neural state. Recent
extensions of the basic Neural Galerkin framework have sought to address
both points. First, the work \cite{BermanPeherstoferSparse} proposes
sparse updates to the parameters at each time step. These updates
can be viewed as projections onto a subspace of the neural tangent
space. One drawback of this approach is that it relies on an incoherence
assumption for the neural Jacobian. Moreover, due to the changing
sparsity pattern of the update, the parameter dynamics cannot be viewed
as smooth in time. Second, the work \cite{WEN2024134129} has considered
coupling the particle and parameter dynamics, with the goal of maintaining
that the particles sample from an evolving sampling measure induced
by the neural state. However, the approach in this work does not define
a coupled system of smooth ODEs for the particles and parameters.

Our work seeks to address both of these computational questions. First,
we propose an alternative approach to sketching the neural evolution
based on a ridge regression perspective \cite{Martinsson_Tropp_2020,murray2023randomizednumericallinearalgebra}.
This reduces the cost of the integrator from $O(pN^{2})$ to $O(pN\log N+pn^{2})$,
where $n$ is a sketch dimension that can be viewed as proportional
(up to log factors) to the numerical rank of neural Jacobian. Importantly,
we can use a fixed oblivious sketch matrix, allowing us to formulate
smooth dynamics for the parameters that we call the minimal neural
evolution (MNE).

To address the coupling of the parameter and collocation point dynamics,
we focus specifically on applications to diffusion processes, where
the score function \cite{song2021scorebased} of the evolving state
suggests an intuitive evolution of collocation points alongside the
network parameters, similar to recent work on Fokker-Planck evolution
\cite{Boffi_2023,BoffiActive} which retrains neural network parameters
at each time step, rather than propagating them directly. By specifying
joint smooth dynamics for the parameters and collocation points, we
can leverage high-order adaptive integration methods for efficient
and accurate propagation of both parameters and collocation points
over time.

Interestingly, we do not encounter any analog of the tangent space
collapse phenomenon \cite{zhang2024sequentialintimetrainingnonlinearparametrizations}
identified as a danger in Neural Galerkin schemes. Relatedly, we believe
that further analysis of the generalization error of the MNE is an
interesting topic for future study.

A key application of our work is to sampling problems, where MNE provides
a method for generating samples from a probability distribution given
its density function but no direct training data. In this setting,
a natural \emph{a priori} specification of the evolution of the collocation
points is possible, further simplifying the integration.

A large body of recent work has focused on the application of ideas
from generative modeling, powered by neural networks, to this sampling
problem. We do not attempt a comprehensive review, but we refer to
the recent work \cite{albergo2025netsnonequilibriumtransportsampler}
for a thorough accounting and comparison of recent methods. We point
out that many recent works such as Liouville flow \cite{LiouvilleFlow},
NETS \cite{albergo2025netsnonequilibriumtransportsampler}, and AFT
/ CRAFT \cite{arbel2021,matthews2022} are based on the identification
of an `annealing' trajectory between the target density and an easier
reference density, induced by linear interpolation of the log-densities.
Others seek to freely optimize over density trajectories \cite{doucet2022scorebased,LangevinDiffusionVI,StochasticNormalizing}.
By contrast, our approach relies on the trajectory which connects
the target density to a Gaussian distribution via Ornstein-Uhlenbeck
diffusion, like diffusion models \cite{pmlr-v37-sohl-dickstein15,NEURIPS2020_Ho,song2021scorebased}.
Relatedly, the recent work \cite{vargas2023denoising} also learns
this trajectory, through alternative techniques.

One key advantage of the diffusion approach is that it allows for
natural generalization to the task of conditional sampling and marginalization,
which are of interest in many scientific applications, and can be
similarly accommodated in generative models based on diffusion \cite{ho2021classifierfree,ho2022video,ho2021cascaded}.
We illustrate our approach to these more general tasks and validate
it with numerical experiments on a high-dimensional target distribution.

Unlike most other approaches to neural sampling, our approach to sampling
relies on a preprocessing step where the target density is fit with
a neural surrogate. Importantly, this preprocessing step does not
rely on having exact samples from the target density. We believe that
this could be an interesting starting point for other approaches as
well. In particular, many target densities in Bayesian inference,
which are of moderate dimension and therefore of particular interest
as potential loci of `neural advantage' in sampling (see, e.g., \cite{karamanis2022accelerating,karamanis2022pocomc,pmlr-v139-dai21a,doi:10.1073/pnas.2309624121,TRENF}),
are characterized by extremely expensive density evaluations. For
such distributions, relying on online density evaluations within a
neural network training subroutine may be completely impractical,
and training a neural surrogate on a limited set of representative
collocation points (computed offline) may be a useful preprocessing
step for other methods as well. (One of our example tasks in this
work, sampling the Bayesian posterior for hyperparameters in Gaussian
process regression \cite{rasmussen_gaussian_2006,pmlr-v118-lalchand20a},
can be viewed as such a problem. For large datasets, each density
evaluation requires the evaluation of an expensive matrix determinant.)

Additionally, we extend our methodology to conditional sampling and
marginalization, and we demonstrate how to systematically remove parametric
approximation biases in the generated samples. We also demonstrate
that the bias due to parametric approximation can be systematically
removed from our samples (either full and conditional), similar to
the recent work \cite{albergo2025netsnonequilibriumtransportsampler},
via particle reweighting similar to diffusion Monte Carlo \cite{becca2017quantum},
with optional resampling to avoid population collapse.

To validate our approach, we present numerical experiments on a variety
of diffusion-driven problems, including Langevin dynamics with nonequilibrium
forcing, fully Bayesian Gaussian process regression, and sampling
from an Allen-Cahn potential. These experiments demonstrate the effectiveness
and scalability of MNE, showing that it can systematically improve
solution accuracy while maintaining computational feasibility in high-dimensional
settings.

This paper is structured as follows. In Section 2, we formalize the
minimal neural evolution (MNE) framework. Section 3 presents the application
of MNE to diffusion processesa and, in turn, to sampling and marginalization.
Section 4 provides numerical results demonstrating the accuracy and
computational efficiency of the method.

\subsection*{Acknowledgments}

This work was supported in part by a Sloan Research Fellowship and
by the U.S. Department of Energy, Office of Science, Office of Advanced
Scientific Computing Research's Applied Mathematics Competitive Portfolios
program under Contract No. AC02-05CH11231.

\section{Minimal neural evolution \label{sec:Minimal-neural-evolution}}

Consider a function evolution dictated by the equation

\begin{equation}
\partial_{t}u_{t}(x)=\mathcal{A}_{t}[u_{t}](x),\quad x\in\R^{d},\label{eq:evol}
\end{equation}
 where $\mathcal{A}_{t}$ is an evolution operator that is possibly
nonlinear and time-dependent.

Given a parametrization $u_{\theta}(x)$ of a function in terms of
parameters $\theta\in\R^{p}$, we want to solve (\ref{eq:evol}) within
the parametric manifold $\{u_{\theta}\,:\,\theta\in\R^{p}\}$, or
equivalently, seek an evolution for the parameters $\theta=\theta(t)$.

\subsection{Preliminary derivation \label{subsec:Preliminary-derivation}}

To determine this evolution, we will insist that (\ref{eq:evol})
is satisfied exactly at a collection of collocation points $X=(x_{1},\ldots,x_{N})$,
or almost exactly where a regularization parameter specifies the softness
of these constraints. The collection $X=X(t)$ may be time-dependent
with either a predetermined evolution or an evolution that is coupled
to that of the parameters $\theta(t)$. We will use the notation $f(X)\in\R^{N}$
to denote the batched evaluation of a function $f:\R^{d}\ra\R$ at
all points in $X$. More generally, if $f(x)$ is vector-valued, then
$f(X)$ will be matrix-valued with $N$ columns.

Now for large overparametrized networks with $p\gg N$, the evolution
of $\theta$ may be underspecified by the constraint at $X$. Therefore,
for a small time increment, we propose to find the smallest change
in the neural network parameters that is compatible with solving (\ref{eq:evol})
at the collocation points.

Specifically, supposing that $\theta=\theta(t)$ is the currrent value
of the parameters at some time $t$ and $X=X(t)$ the current position
of the collocation points, we seek to find the rate of change $\dot{\theta}=\dot{\theta}(t)$
for the parameters that satisfies: 
\begin{equation}
\underset{\dot{\theta}\in\R^{p}}{\text{minimize}}\left\Vert \left[\nabla_{\theta}u_{\theta}(X)\right]^{\top}\dot{\theta}-\mathcal{A}_{t}[u_{\theta}](X)\right\Vert ^{2}+\lambda\Vert\dot{\theta}\Vert^{2},\label{eq:ridge}
\end{equation}
 where $\Vert\,\cdot\,\Vert$ denotes the usual 2-norm and $\lambda>0$
is a regularization parameter. Note that to derive this optimization
problem we have used the chain rule $\frac{d}{dt}\left[u_{\theta(t)}(x)\right]=\left\langle \nabla_{\theta}u_{\theta(t)}(x),\dot{\theta}(t)\right\rangle $.
For shorthand we will write 
\begin{equation}
\Phi_{\theta}(x)=\nabla_{\theta}u_{\theta}(x),\label{eq:featuremap}
\end{equation}
 which can be viewed as a feature map in the sense of the neural tangent
kernel (NTK) \cite{NTK}.

Now the standard linear-algebraic manipulations for the ridge regression
problem (\ref{eq:ridge}) reveal that the optimizer $\dot{\theta}$
is given by 
\begin{equation}
\dot{\theta}=\Phi_{\theta}(X)\left[\Phi_{\theta}(X)^{\top}\Phi_{\theta}(X)+\lambda\mathbf{I}_{p}\right]^{-1}\mathcal{A}_{t}[u_{\theta}](X).\label{eq:MNE}
\end{equation}
 We can view (\ref{eq:MNE}) as an ODE specifying the evolution of
$\theta(t)$. The evolution of $X(t)$ remains to be specified and
the appropriate choice will be context-dependent, and we will discuss
the choice below. Modulo this specification, as mentioned in the introduction,
(\ref{eq:MNE}) can be viewed as a Neural Galerkin scheme \cite{Bruna_Peherstorfer_Vanden-Eijnden_2024}
in which the same regularization parameter $\lambda$ has been used
to shift the metric. 

Note that by comparing terms in (\ref{eq:ridge}), it makes sense
to take 
\begin{equation}
\lambda=N\ve^{2},\label{eq:regparam}
\end{equation}
 where $\ve>0$ roughly specifies an order of magnitude for the pointwise
error tolerance.

Also observe that the matrix 
\[
K_{\theta}(X,X):=\Phi_{\theta}(X)^{\top}\Phi_{\theta}(X)
\]
 appearing in (\ref{eq:MNE}) is the kernel matrix of the NTK \cite{NTK}
evaluated on our collocation points. Modulo regularization, (\ref{eq:MNE})
seeks the update to the neural network parameters which is smallest
in the sense of the NTK inner product while maintaining the evolution
constraints at the collocation points. Note that this interpretation
is not dependent on the evolution staying within some kernel learning
regime \cite{NTK}.

\subsection{Sketching \label{subsec:Sketching}}

The matrix inversion required to implement (\ref{eq:MNE}) is $N\times N$,
by contrast with the $p\times p$ matrix inversion that appears in
the neural Galerkin method \cite{Bruna_Peherstorfer_Vanden-Eijnden_2024}---though,
as mentioned above, this may alternatively be viewed as a least squares
problem. In order to resolve the evolution with high fidelity throughout
the spatial domain, it is of interest to enlarge the size $N$ of
the collocation set $X$ as much as possible.

At the same time, enlarging this collocation set may introduce a significant
amount of redundancy into the constraints due to rank deficiency of
the NTK. This redundancy can be removed efficiently via a sketching
procedure. Specifically, we will let $\Omega\in\R^{n\times N}$ denote
a sketch matrix \cite{Martinsson_Tropp_2020} which reduces the dimension
of a vector from $N$ to $n$ when multiplied from the left. Even
more concretely, we consider as our sketch the subsampled randomized
discrete Hartley transform, which is a partial isometry and is considered
to be the best fast sketch for real matrices \cite{Martinsson_Tropp_2020}:
\[
\Omega v=\frac{1}{\sqrt{n}}\left[\Re\left(\mathbf{F}_{N}[\sigma\odot v]\right)-\Im(\mathbf{F}_{N}[\sigma\odot v])\right]_{\mathcal{I}}.
\]
 Here $\sigma\in\{\pm1\}^{N}$ is a fixed vector of i.i.d. random
signs, $\mathbf{F}_{N}$ denotes the $N$-dimensional discrete Fourier
transform, and $\mathcal{I}$ is a fixed uniformly chosen random subset
of $\{1,\ldots,N\}$ of size $n$.

Then we can consider in place of (\ref{eq:ridge}) the sketched ridge
regression problem 

\begin{equation}
\underset{\dot{\theta}\in\R^{p}}{\text{minimize}}\left\Vert \Omega\left[\nabla_{\theta}u_{\theta}(X)\right]^{\top}\dot{\theta}-\Omega\mathcal{A}_{t}[u_{\theta}](X)\right\Vert ^{2}+\lambda\Vert\dot{\theta}\Vert^{2},\label{eq:sketchedridge}
\end{equation}
 which in turn yields 
\begin{equation}
\dot{\theta}=\tilde{\Phi}_{\theta}(X)\left[\tilde{\Phi}_{\theta}(X)^{\top}\tilde{\Phi}_{\theta}(X)+\lambda\mathbf{I}_{p}\right]^{-1}\Omega\,\mathcal{A}_{t}[u_{\theta}](X),\label{eq:sketchedMNE}
\end{equation}
 where 
\[
\tilde{\Phi}_{\theta}(X):=\Phi_{\theta}(X)\,\Omega^{\top}.
\]
 Notably the cost of forming $\tilde{\Phi}_{\theta}(X)$ is only $O(pN\log N)$
and the total cost of implementing the right-hand side of the sketched
MNE (\ref{eq:sketchedMNE}) is only 
\[
O(pN\log N+pn^{2}+n^{3}),
\]
 where we recall that $N$ is the number of collocation points, $n$
is the sketch dimension, and $p$ is the number of parameters. Typically
a neural network ansatz will be overparametrized in the sense that
the numerical rank of the NTK will be far smaller than $p$. Hence
we can typically take $n<p$ (or even $n\ll p$) to be roughly on
the order of the numerical rank of the NTK \cite{Martinsson_Tropp_2020,murray2023randomizednumericallinearalgebra}.

We will fix a sketch matrix $\Omega$ throughout the dynamics, allowing
us to view the right-hand side of (\ref{eq:sketchedMNE}) as a smooth,
deterministic function of $\theta$ (and possibly $X$), so that standard
solvers can be applied. In particular, we implement the ODE solver
in JAX \cite{jax2018github} using the Diffrax library \cite{kidger2021on}.
We choose Tsitouras' 5/4 method \cite{10.1016/j.camwa.2011.06.002}
(i.e., $\texttt{diffrax.Tsit5}$) as our solver, which is an explicit
Runge-Kutta method with adaptive time-stepping.

Subject to the specification of the dynamics of $X(t)$, we call (\ref{eq:sketchedMNE})
the \emph{minimal neural evolution }(\emph{MNE}).

\section{Application to diffusion processes \label{sec:Application-to-diffusion}}

We will be interested in the density evolution induced by the stochastic
differential equation \cite{Risken_1996}
\begin{equation}
dX_{t}=b_{t}(X_{t})\,dt+\sigma\,\chi_{\mathcal{S}}\odot dB_{t}.\label{eq:sde}
\end{equation}
 Here $b_{t}:\R^{d}\ra\R^{d}$ is a time-dependent drift field, $\sigma>0$
is a diffusion coefficient, and $\chi_{\mathcal{S}}=\sum_{i\in\mathcal{S}}e_{i}$
is a `mask vector' which indicates the subset $\mathcal{S}\subset\{1,\ldots,d\}$
of variables that are subject to diffusion. More general diffusion
terms could be considered, but we restrict ourself to this case for
simplicity since it covers all applications of interest below. Note
with caution that throughout we will use the standard notation $X_{t}$
for the SDE variable, but it should not be confused with our notation
$X=X(t)$ for the collocation points.

\subsection{General discussion \label{subsec:General-discussion}}

First we explain how the density evolution induced by (\ref{eq:sde})
can be approached within the MNE framework.

To begin, observe that the Fokker-Planck equation \cite{Risken_1996}
corresponding to (\ref{eq:sde}) is given by 
\begin{equation}
\partial_{t}\rho_{t}=-\nabla\cdot(\rho_{t}b_{t})+\frac{\sigma^{2}}{2}\Delta_{\mathcal{S}}\rho_{t},\label{eq:fokkerplanck}
\end{equation}
 where $\Delta_{\mathcal{S}}=\sum_{i\in\mathcal{S}}\frac{\partial^{2}}{\partial x_{i}^{2}}$
is a masked Laplacian.

We will not parametrize the evolution of $\rho_{t}$, but rather that
of the energy function $u_{t}=-\log\rho_{t}$. By substitution and
repeated application of the chain rule, one deduces the following
evolution equation for $u_{t}$: 
\begin{equation}
\partial_{t}u_{t}=\nabla\cdot b_{t}-b_{t}\cdot\nabla u_{t}+\frac{\sigma^{2}}{2}\Delta_{\mathcal{S}}u_{t}-\frac{\sigma^{2}}{2}\vert\nabla_{\mathcal{S}}u_{t}\vert^{2},\label{eq:energyevol}
\end{equation}
 where $\nabla_{\mathcal{S}}:=\left(\frac{\partial}{\partial x_{i}}\right)_{i\in\mathcal{S}}$
denotes the masked gradient.

Then we can directly apply MNE to the operator 
\begin{equation}
\mathcal{A}_{t}[u]=\nabla\cdot b_{t}-b_{t}\cdot\nabla u+\frac{\sigma^{2}}{2}\Delta_{\mathcal{S}}u-\frac{\sigma^{2}}{2}\vert\nabla_{\mathcal{S}}u\vert^{2}.\label{eq:mnediffusion}
\end{equation}

We would like to insist that the collocation points $X(t)$ evolve
with the density. To intialize them, we might choose $X(0)$ to be
samples from $\rho_{0}$. However, the consistency of the scheme does
not require us to do so exactly. Either way, it is natural to demand
that $X(t)$ evolve according to a deterministic flow that implements
the density evolution.

Since the Fokker-Planck equation can be viewed as a transport equation
\[
\partial_{t}\rho_{t}=-\nabla\cdot\left(\rho_{t}\left[b_{t}+\frac{\sigma^{2}}{2}\chi_{\mathcal{S}}\odot\nabla u_{t}\right]\right)
\]
 with drift field $b_{t}+\frac{\sigma^{2}}{2}\chi_{\mathcal{S}}\odot\nabla u_{t}$,
we can impose the following evolution for the collocation points:
\[
\dot{X}=b_{t}(X)+\frac{\sigma^{2}}{2}\chi_{\mathcal{S}}\odot\nabla u_{\theta}(X),
\]
 where the entrywise product is broadcasted suitably. Together with
equation (\ref{eq:sketchedMNE}) for $\dot{\theta}$, this completes
the specification of a system of ODEs for $(\theta(t),X(t))$.

However, in the application of MNE to marginalization and sampling
considered below, we will consider an alternative (\emph{a priori})
specification of $X(t)$ which need not be solved for jointly with
$\theta(t)$.

\subsection{Sampling and marginalization \label{subsec:Sampling-and-marginalization}}

We can apply the density evolution framework of the preceding Section
\ref{subsec:General-discussion} to the problem of sampling a probability
density 
\[
\rho(x)\propto e^{-u(x)},
\]
 where $u(x)$ is known.

Our approach also allows for marginal estimation and conditional sampling.
Specifically, for any subset of variables $\mathcal{S}\subset\{1,\ldots,d\}$,
let $\mathcal{S}'$ denote its complement, and identify $x=(x_{\mathcal{S}},x_{\mathcal{S}'})$
in the notation, assuming for simplicity that the variables have been
so ordered. Then let 
\[
\rho_{\mathcal{S}'}(x_{\mathcal{S}'}):=\int_{\R^{\vert\mathcal{S}\vert}}\rho(x)\,dx_{\mathcal{S}}
\]
 denote the marginal density. We will show how to construct an approximation
of the marginal energy function $u_{\mathcal{S}'}$ which is defined
up to a constant shift by 
\[
\rho_{\mathcal{S}'}(x_{\mathcal{S}'})\propto e^{-u_{\mathcal{S}'}(x_{\mathcal{S}'})}.
\]

Moreover let 
\[
\rho(x_{\mathcal{S}}\,\vert\,x_{\mathcal{S}'})=\frac{\rho(x)}{\rho_{\mathcal{S}'}(x_{\mathcal{S}'})}
\]
 denote the conditional density for the variables in $\mathcal{S}$
given those in $\mathcal{S}'$. We will show how to draw samples $x_{\mathcal{S}}$
from this distribution for arbitrary given $x_{\mathcal{S}'}$. In
particular, we will show how the bias due to our parametric approximations
can be systematically removed from these samples. Note that in the
special case $\mathcal{S}=\{1,\ldots,d\}$, we recover the capacity
to draw samples from the full target distribution $\rho$.

\subsubsection{Preprocessing \label{subsec:Preprocessing}}

As a first step, we obtain some initial collocation points $X_{\mathrm{init}}=(x_{1},\ldots,x_{N})$
adapted to the target density $\rho$. Importantly, we do not need
to insist that $X_{\mathrm{init}}$ consists of samples from $\rho$.
Instead, we can always obtain $X_{\mathrm{init}}$ by running a few
iterations of a Markov chain Monte Carlo (MCMC) sampler with respect
to the target density, so that $X_{\mathrm{init}}$ attains coverage
of the `typical set' \cite{Betancourt_2018} of $\rho$. For example,
even if multimodality prevents rapid mixing of MCMC, it is not necessarily
a problem if the initial samples $X_{\mathrm{init}}$ exhibit an improper
balance between modes, as long as all modes are represented.

The next step is to fit $u$ with a neural network surrogate $u_{\theta_{\mathrm{init}}}$
on our initial collocation points. We achieve this by minimizing the
error of the score function on our collocation points 
\begin{equation}
\mathcal{L}(\theta)=\frac{1}{N}\sum_{i=1}^{N}\vert\nabla u_{\theta}(x_{i})-\nabla u(x_{i})\vert^{2},\label{eq:scorematching}
\end{equation}
 i.e., setting $\theta_{\mathrm{init}}=\underset{\theta\in\R^{p}}{\text{argmin}}\ \mathcal{L}(\theta)$.

\subsubsection{MNE specification \label{subsec:MNE-specification}}

In the context of the general presentation in Section \ref{subsec:General-discussion}
we make the choice 
\begin{equation}
b_{t}(x)=-\gamma\,\chi_{\mathcal{S}}\odot x\label{eq:drift}
\end{equation}
 for the drift field, where $\gamma>0$ is a damping parameter. Then
the operator $\mathcal{A}=\mathcal{A}_{t}$ specified by (\ref{eq:mnediffusion})
is time-independent.

For this choice of $b_{t}$, the density $\rho_{t}$ for the diffusion
process (\ref{eq:sde})-(\ref{eq:fokkerplanck}) converges rapidly
to the product density 
\[
\rho_{\infty}(x)\propto e^{-\frac{1}{\sigma^{2}/\gamma}\Vert x_{\mathcal{S}}\Vert^{2}}\,\rho_{\mathcal{S}'}(x_{\mathcal{S}'}),
\]
 i.e., the product of the marginal density $\rho_{\mathcal{S}'}$
for the variables in $\mathcal{S}'$ with a normal distribution $\mathcal{N}\left(0,\frac{\sigma^{2}}{2\gamma}\mathbf{I}_{\vert\mathcal{S}\vert}\right)$
for the variables in $\mathcal{S}$.

Since the dynamics become smoother at larger times, we introduce the
time change $t(s)=\frac{1}{2}s^{2}$ and evolve with respect to the
transformed variable $s$, defining $\tilde{u}_{s}=u_{t(s)}$. Then
by the chain rule, together with (\ref{eq:energyevol})-(\ref{eq:mnediffusion}),
we have that 
\begin{equation}
\partial_{s}\tilde{u}_{s}(x)=s\mathcal{A}[\tilde{u}_{s}](x),\label{eq:timechange}
\end{equation}
 and we can apply the MNE induced by the choice of operator $\tilde{\mathcal{A}_{s}}[u]=s\mathcal{A}[u]$
to implement the evolution of $\tilde{u}_{s}$ in the transformed
time variable $s$, deducing an evolution of parameters $\theta(s)$
such that $\tilde{u}_{s}\approx u_{\theta(s)}$. We choose $\theta(0)=\theta_{\mathrm{init}}$
as our initial condition, where $\theta_{\mathrm{init}}$ consists
of the parameters of our surrogate model as computed in Section \ref{subsec:Preprocessing}.

Then we view $u_{0}=u_{\theta(0)}$ as the exact initial condition
for our dynamics (\ref{eq:timechange}). We will let 
\[
\rho_{s}(x)\propto e^{-\tilde{u}_{s}(x)}
\]
 denote the density induced by $\tilde{u}_{s}$ over the course of
the exact evolution (\ref{eq:timechange}), and we will overload the
notation $\rho_{s}$ to indicate corresponding conditional distributions
as well, as shall be clear from context.

\subsubsection{Collocation points \label{subsec:Collocation-points}}

With initial collocation points $X_{\mathrm{init}}$ obtained as specified
in Section \ref{subsec:Preprocessing}, we obtain a trajectory $X(t)$
via the formula

\begin{equation}
X(t)=(\mathbf{1}-\chi_{\mathcal{S}})\odot X_{\mathrm{init}}+\chi_{\mathcal{S}}\odot\left[e^{-\gamma t}X_{\mathrm{init}}+\sqrt{\frac{\sigma^{2}}{2\gamma}\left(1-e^{-2\gamma t}\right)}\ Z\right],\label{eq:diffusiontraj}
\end{equation}
 where the columns of $Z=(z_{1},\ldots,z_{N})$ are i.i.d. samples
from $\mathcal{N}(0,\mathbf{I}_{d})$. Note that only the variables
in $\mathcal{S}$ are changing over this trajectory.

The trajectory (\ref{eq:diffusiontraj}) is chosen so that if the
columns of $X_{\mathrm{init}}$ are drawn from some initial sampling
density $\rho_{\mathrm{init}}$ (which need not coincide with the
target density $\rho$), then the columns of $X(t)$ sample exactly
from the density obtained by propagating $\rho_{\mathrm{init}}$ according
to the Fokker-Planck evolution specified by our drift (\ref{eq:drift})
up to time $t$. Importantly, the dependence of $X(t)$ on $t$ is
smooth, and by contrast it is undesirable to produce a trajectory
by solving the SDE (\ref{eq:sde}). The same type of sample trajectory
is commonly used for training score-based diffusion models \cite{song2021scorebased}.

More properly, after the time change to $s$ via $t(s)=\frac{1}{2}s^{2}$,
we choose $\tilde{X}(s)=X(t(s))$ as our set of evolving collocation
points.

\subsubsection{Marginalization \label{subsec:Marginalization}}

Note that if $\vert\mathcal{S}\vert<n$, then we can deduce the energy
function $u_{\mathcal{S}'}$ of the marginal distribution $\rho_{\mathcal{S}'}\propto e^{-u_{\mathcal{S}'}}$
via the evaluation
\[
u_{\mathcal{S}'}(x_{\mathcal{S}'}):=u_{\infty}(0,x_{\mathcal{S}'}).
\]
 In practice we approximate $u_{\infty}\approx u_{\theta(s_{\mathrm{f}})}$,
where $s_{\mathrm{f}}$ is the final time of our evolution.

\subsubsection{Sampling \label{subsec:Sampling}}

We can use the learned MNE to draw samples $x_{\mathcal{S}}$ from
the conditional distribution $\rho(\,\cdot\,\vert\,x_{\mathcal{S}'})$,
for arbitrary given $x_{\mathcal{S}'}$. In the case where $\mathcal{S}=\{1,\ldots,d\}$,
we simply draw samples from the full distribution $\rho$.

To achieve this, consider the reverse diffusion \cite{song2021scorebased}:
\begin{equation}
dX_{s}=(s_{\mathrm{f}}-s)\,\chi_{\mathcal{S}}\odot\left[\gamma X_{s}-\sigma^{2}\nabla\tilde{u}_{s_{\mathrm{f}}-s}(X_{s})\right]\,ds+(s_{\mathrm{f}}-s)^{1/2}\,\sigma\,\chi_{\mathcal{S}}\odot dB_{s},\quad s\in[0,s_{\mathrm{f}}].\label{eq:rev}
\end{equation}
 Again note with caution that we use the standard notation $X_{t}$
for the SDE variable, but it should not be confused with our notation
$X=X(t)$ for the collocation points.

The Fokker-Planck evolution of (\ref{eq:rev}) reverses the Fokker-Planck
evolution (\ref{eq:fokkerplanck}) of the forward diffusion (\ref{eq:sde}),
taking into account the time change $t=t(s)$, and therefore if $X_{0}\sim\rho_{s_{\mathrm{f}}}$
then $X_{\mathrm{s}_{f}}\sim\rho_{0}$.

Since these dynamics leave the variables in $\mathcal{S}'$ unaltered,
we can alternatively consider the simpler dynamics for the remaining
variables $Y_{s}=[X_{s}]_{\mathcal{S}}$: 
\begin{equation}
dY_{s}=(s_{\mathrm{f}}-s)\left[\gamma Y_{s}-\sigma^{2}\nabla_{\mathcal{S}}\tilde{u}_{s_{\mathrm{f}}-s}(Y_{s},x_{\mathcal{S}'})\right]\,ds+(s_{\mathrm{f}}-s)^{1/2}\,\sigma\,dB_{s},\quad s\in[0,s_{\mathrm{f}}],\label{eq:revmask}
\end{equation}
 where by some abuse of notation $B_{s}$ here indicates a Brownian
motion of dimension $\vert\mathcal{S}\vert$. This diffusion must
transport the conditional distribution $\rho_{s_{\mathrm{f}}}(\,\cdot\,\vert\,x_{\mathcal{S}'})$
to $\rho_{0}(\,\cdot\,\vert\,x_{\mathcal{S}'})$ for any given $x_{\mathcal{S}'}$.

Under the assumption that $s_{\mathrm{f}}$ is sufficiently large,
we have that the law of the conditional distribution $\rho_{s_{\mathrm{f}}}(\,\cdot\,\vert\,x_{\mathcal{S}'})$
is approximately $\mathcal{N}\left(0,\frac{\sigma^{2}}{2\gamma}\mathbf{I}_{\vert\mathcal{S}\vert}\right)$,
independent of $x_{\mathcal{S}'}$. Moreover, $\rho_{0}(\,\cdot\,\vert\,x_{\mathcal{S}'})\approx\rho(\,\cdot\,\vert\,x_{\mathcal{S}'})$,
provided that our initial condition $\tilde{u}_{0}=u_{\theta(0)}$
for the MNE approximates the target energy function $u$.

Therefore, to draw an approximate sample $y$ from the conditional
distribution $\rho(\,\cdot\,\vert\,x_{1:d'})$, we would like to draw
$Y_{0}\sim\mathcal{N}\left(0,\frac{\sigma^{2}}{2\gamma}\mathbf{I}_{\vert\mathcal{S}\vert}\right)$
and then solve the SDE (\ref{eq:revmask}) to obtain our sample as
$Y_{s_{\mathrm{f}}}$. However, since we do not have access to the
exact evolution $\tilde{u}_{s}$ but rather only our neural approximation
$u_{\theta(s)}$, we instead obtain our sample by solving, obtained
from (\ref{eq:revmask}) via this substitution: 
\begin{equation}
dY_{s}=(s_{\mathrm{f}}-s)\left[\gamma Y_{s}-\sigma^{2}\nabla_{\mathcal{S}}u_{\theta(s_{\mathrm{f}}-s)}(Y_{s},x_{\mathcal{S}'})\right]\,ds+(s_{\mathrm{f}}-s)^{1/2}\,\sigma\,dB_{s},\quad s\in[0,s_{\mathrm{f}}].\label{eq:revmask-1}
\end{equation}

\subsubsection{Unbiasing \label{subsec:Unbiasing}}

There are several sources of bias in the sampling procedure described
above due to the parametric approximation error both of the initial
condition $u_{\theta(0)}\approx u$ and the evolution itself. All
of these can be systematically removed.

Note that if the MNE dynamics $u_{\theta(s)}$ were to exactly satisfy
(\ref{eq:timechange}), and if $Y_{0}$ was drawn exactly from $\rho_{s_{\mathrm{f}}}(\,\cdot\,\vert\,x_{\mathcal{S}'})$,
then $Y_{s}$ computed via (\ref{eq:revmask-1}) would be an exact
sample from $\rho_{s_{\mathrm{f}}-s}(\,\cdot\,\vert\,x_{\mathcal{S}'})$
for all $s\in[0,s_{\mathrm{f}}]$. More generally, we can use the
discrepancy by which the equality (\ref{eq:timechange}) fails to
hold to evolve weights for our samples. The same idea has also appeared
in the recent work \cite{albergo2025netsnonequilibriumtransportsampler}.

To wit, we couple the SDE (\ref{eq:revmask}) for $Y_{s}$ to an evolution
for the log-weight of the sample 
\begin{equation}
\dot{w}(s)=(s_{\mathrm{f}}-s)\mathcal{A}[u_{\theta(s_{\mathrm{f}}-s)}](Y_{s},x_{\mathcal{S}'})+\frac{d}{ds}u_{\theta(s_{\mathrm{f}}-s)}(Y_{s},x_{\mathcal{S}'}),\label{eq:logweightdyn}
\end{equation}
 where $w(0)=0.$ In so doing we can obtain a weighted sample $(y,w)=(Y_{s_{\mathrm{f}}},w(s_{\mathrm{f}}))$.
Repeating this procedure $M$ times in parallel we can obtain a collection
of weighted samples $(y^{(i)},w^{(i)})_{i=1}^{M}$, from which empirical
estimates according to the conditional distribution $\rho_{0}(\,\cdot\,\vert\,x_{\mathcal{S}'})$
can be estimated using the weighted empirical measure 
\begin{equation}
\frac{\sum_{i=1}^{M}e^{w^{(i)}}\delta_{y^{(i)}}}{\sum_{i=1}^{M}e^{w^{(i)}}}.\label{eq:weightedempirical}
\end{equation}
 During the propagation, if the effective sample size \cite{ESS}
of the weighted ensemble drops below $\alpha M$ for some parameter
$\alpha\in(0,1)$, we can perform unbiased resampling to avoid population
collapse. Although we have implemented and validated this capacity,
we never need to rely on it in our experiments reported below. In
principle, MCMC steps with respect to the target $e^{-u_{\theta(s_{\mathrm{f}}-s)}}$
could also be included after resampling, as is done in several recent
works \cite{arbel2021,matthews2022}.

Finally, we must address the bias due to the facts that (1) the law
of $\rho_{s_{\mathrm{f}}}(\,\cdot\,\vert\,x_{\mathcal{S}'})$ is not
exactly $\mathcal{N}\left(0,\frac{\sigma^{2}}{2\gamma}\mathbf{I}_{\vert\mathcal{S}\vert}\right)$,
i.e., $Y_{0}$ is not exactly drawn from $\rho_{s_{\mathrm{f}}}(\,\cdot\,\vert\,x_{\mathcal{S}'})$,
and (2) the law of $\rho_{0}(\,\cdot\,\vert\,x_{\mathcal{S}'})$ is
not exactly $\rho(\,\cdot\,\vert\,x_{\mathcal{S}'})$. We can address
both issues by importance sampling.

To wit, after drawing $Y_{0}\sim\mathcal{N}\left(0,\frac{\sigma^{2}}{2\gamma}\mathbf{I}_{\vert\mathcal{S}\vert}\right)$,
we can initialize 
\[
w(0)=\frac{\gamma}{\sigma^{2}}\Vert Y_{0}\Vert^{2}-u_{s_{\mathrm{f}}}(Y_{0},x_{\mathcal{S}'})
\]
 for joint dynamics of $(Y_{s},w(s))$ specified by (\ref{eq:revmask-1})
and (\ref{eq:logweightdyn}). Then after solving these dynamics up
to time $s_{\mathrm{f}}$, we endow our sample $y=Y_{s_{\mathrm{f}}}$
with the modified weight 
\[
w=w(s_{\mathrm{f}})+u_{s_{\mathrm{f}}}(Y_{s_{\mathrm{f}}},x_{\mathcal{S}'})-u(Y_{s_{\mathrm{f}}},x_{\mathcal{S}'}).
\]
 Again, we may repeat this procedure in parallel to draw many weighted
samples $(y^{(i)},w^{(i)})_{i=1}^{M}$ and estimate expectations using
the weighted empirical measure (\ref{eq:weightedempirical}).

\section{Numerical experiments \label{sec:Numerical-experiments}}

Now we describe several numerical experiments validating the MNE framework
for Fokker-Planck evolution and conditional sampling. Our implementation
is in JAX \cite{jax2018github} with double precision, and to solve
ODEs we use the Diffrax library \cite{kidger2021on}. As mentioned
above we choose Tsitouras' 5/4 method \cite{10.1016/j.camwa.2011.06.002}
(i.e., $\texttt{diffrax.Tsit5}$) as our solver, which is an explicit
Runge-Kutta method with adaptive time-stepping, using relative and
absolute tolerances of $10^{-3}$ and $10^{-6}$, respectively. All
experiments were run on a single A100 GPU. The cost of running the
ODE solver for the MNE is the main computational cost of our methodology.
Each such solve needed in Section \ref{subsec:Langevin-dynamics-with}
below was completed in about 30 seconds or below. Each such solve
needed in Sections \ref{subsec:Bayesian-inverse-problem} and \ref{subsec:Allen-Cahn-potential}
was completed in about 10 seconds or below.

\subsection{Langevin dynamics with nonequilibrium forcing \label{subsec:Langevin-dynamics-with}}

We consider the diffusion induced by the following SDE

\[
\begin{cases}
dQ_{t}=P_{t}\,dt\\
dP_{t}=F_{t}(Q_{t})\,dt\,+\sigma\,dB_{t},
\end{cases}
\]
 which are the underdamped Langevin dynamics for a particle with unit
mass, time-dependent force term $q\mapsto F_{t}(q)$, and diffusion
coefficient $\sigma$. The standard form (\ref{eq:sde}) can be recovered
by viewing $X_{t}=(Q_{t},P_{t})$ and $x=(q,p)$. 

As an illustrative example for the MNE we consider the case of $x\in\R^{1+1}$
with forcing specified by 

\[
F_{t}(q)=-q+e^{-q^{2}/2}\cos(t)
\]
 and $\sigma=0.1$.

The Fokker-Planck equation (\ref{eq:fokkerplanck}) is solved using
the MNE approach of Section \ref{sec:Minimal-neural-evolution}, including
the strategy for evolving the collocation point positions $X(t)$.
For our parametrization $\theta\mapsto u_{\theta}$ we use a MLP architecture
with two hidden layers of width 128 and cosine activations.

We take the Gaussian initial condition 
\[
u_{0}=\frac{1}{2}\left[(q-1)^{2}+p^{2}\right],
\]
 and to construct $X(0)$ we simply draw $N=10^{4}$ samples from
$\rho_{0}\propto e^{-u_{0}}$.

To initialize the MNE we must find parameters $\theta(0)$ such that
$u_{\theta(0)}\approx u_{0}$. We achieve this by optimizing the score-matching
objective (\ref{eq:scorematching}) using the Adam optimizer \cite{adam}
with default parameters and learning rate $10^{-2}$. For this example,
since it is possible to fit to a high accuracy, the initial fit error
is ignored and we simply treat $u_{\theta(0)}$ as the exact initial
condition of the dynamics (\ref{eq:energyevol}).

\begin{figure}
\begin{centering}
\includegraphics[scale=0.33]{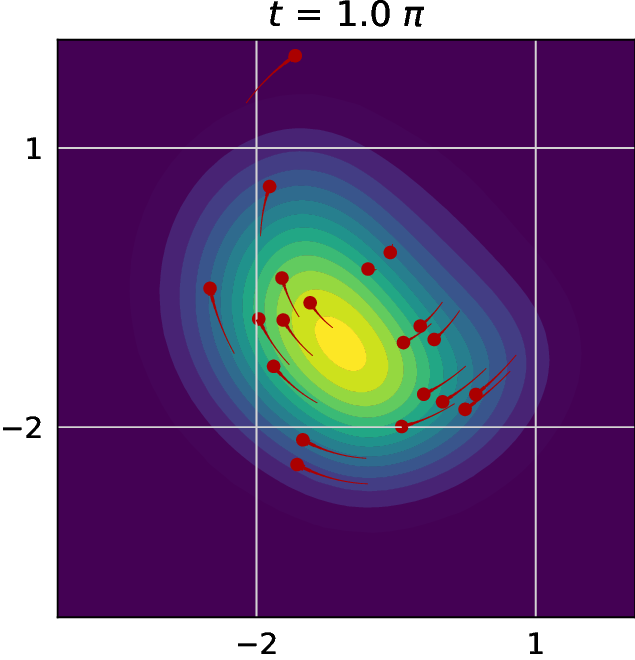} $\ $ \includegraphics[scale=0.33]{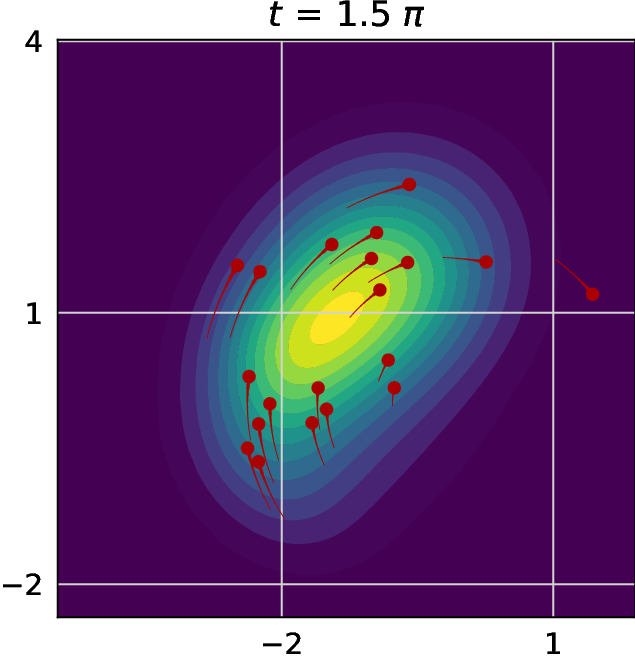}
$\ $ \includegraphics[scale=0.33]{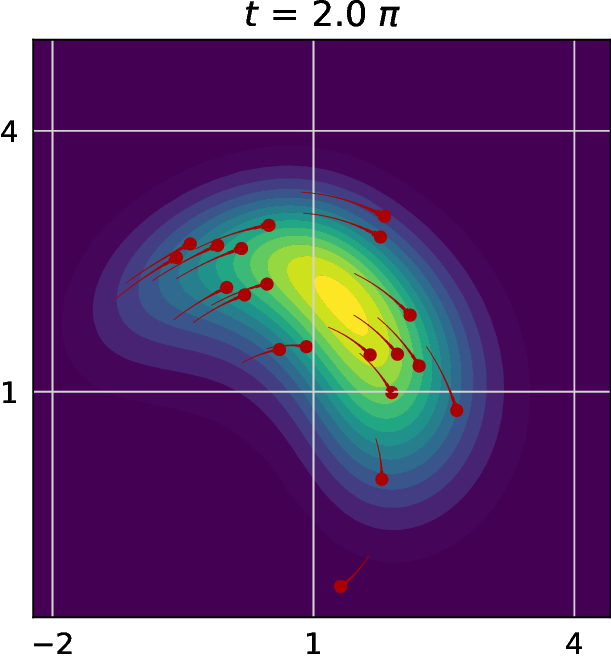} $\ $ \includegraphics[scale=0.33]{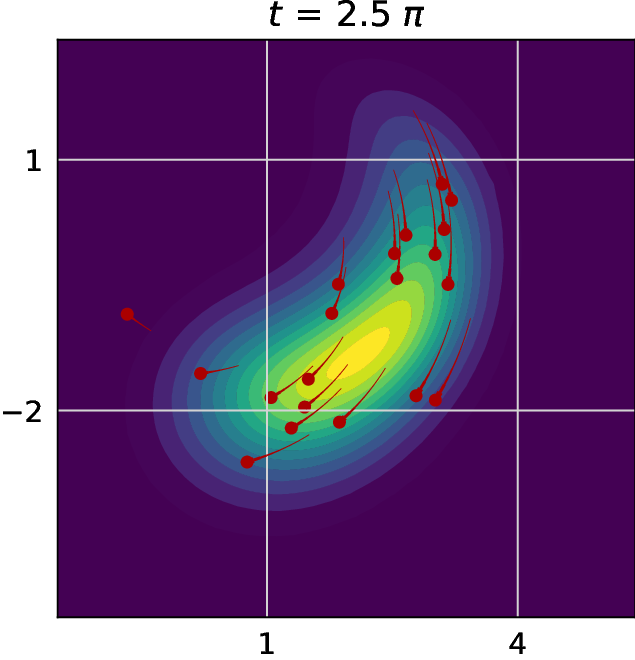}\vspace{3mm}
\par\end{centering}
\begin{centering}
\includegraphics[scale=0.33]{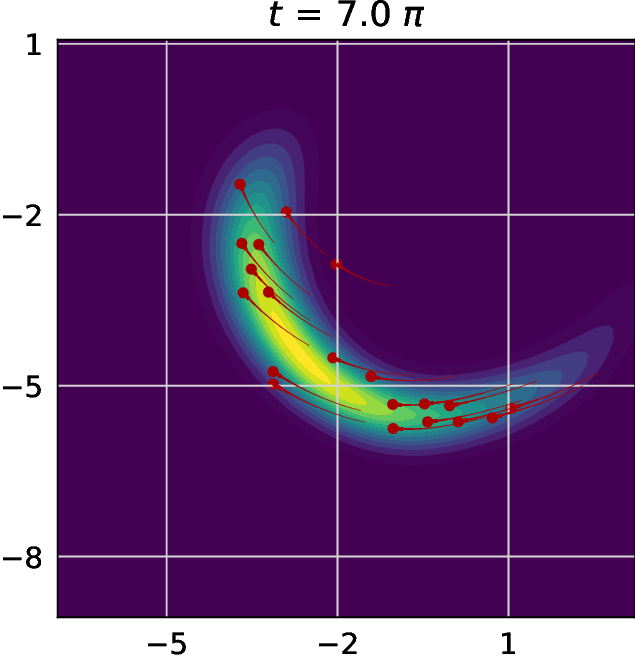} $\ $ \includegraphics[scale=0.33]{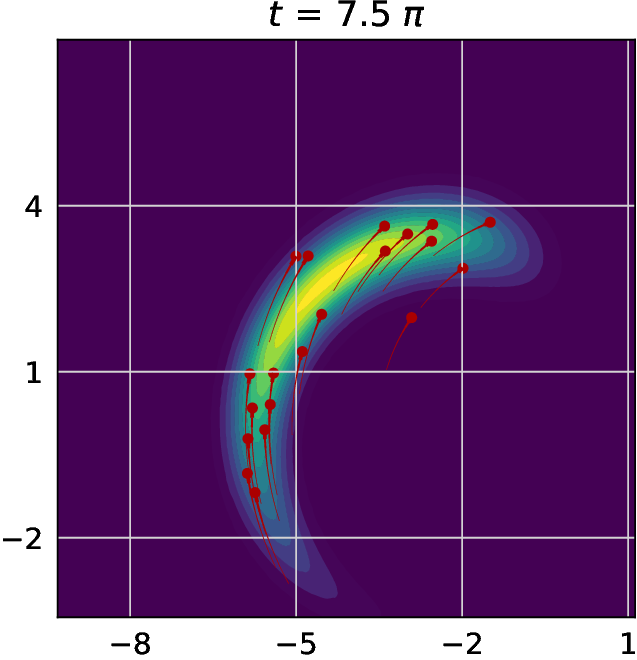}
$\ $ \includegraphics[scale=0.33]{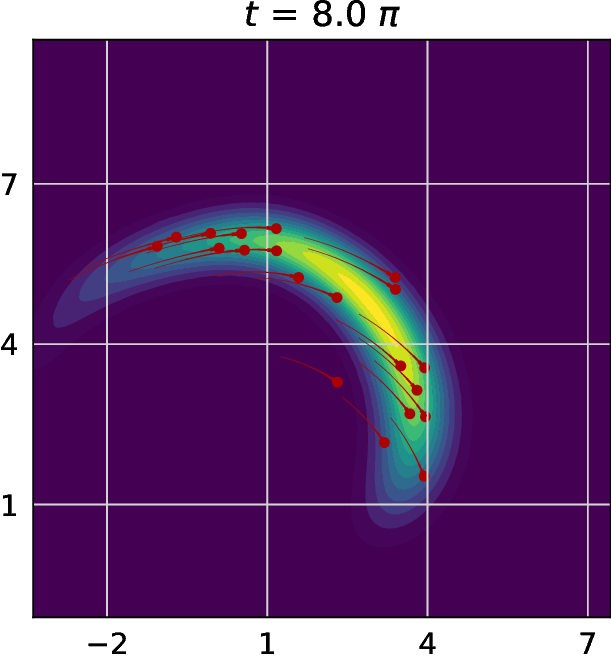} $\ $ \includegraphics[scale=0.33]{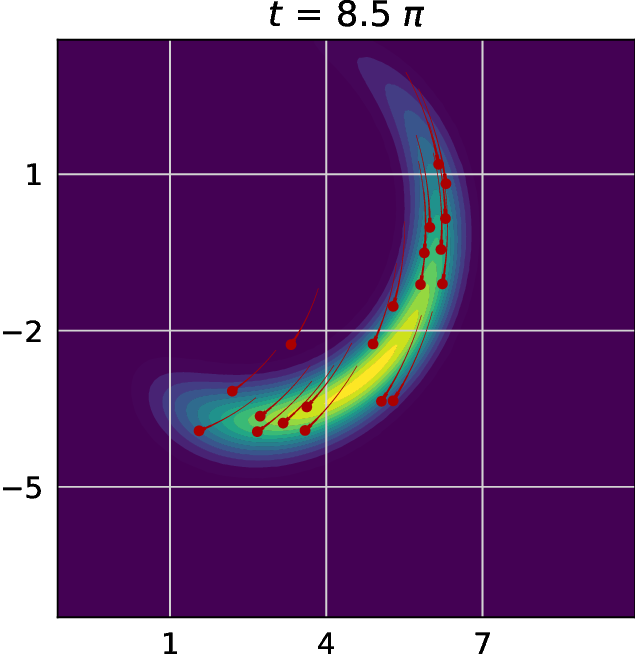}
\par\end{centering}
\caption{RMS residual (\ref{eq:rmseres}) over the MNE trajectory for the example
described in Section \ref{subsec:Langevin-dynamics-with}, plotted
for various choices of the tolerance $\protect\ve$ (cf. (\ref{eq:regparam})).
The sketch dimension is $n=700$ for each curve. \label{fig:movie}}
\end{figure}

Illustrative snapshots of the density evolution and particle trajectories
$X(t)$ are shown in Figure \ref{fig:movie} for several times in
the range $t\in[0,30]$. These were obtained with regularization parameter
$\ve=10^{-6}$ (cf. (\ref{eq:regparam})) and sketch dimension $n=800$.
We will validate the accuracy for these choices below. Animated GIFs
of the evolution are available at 
\[
\texttt{quantumtative.github.io/mne}
\]

\begin{figure}

\begin{centering}
\includegraphics[scale=0.5]{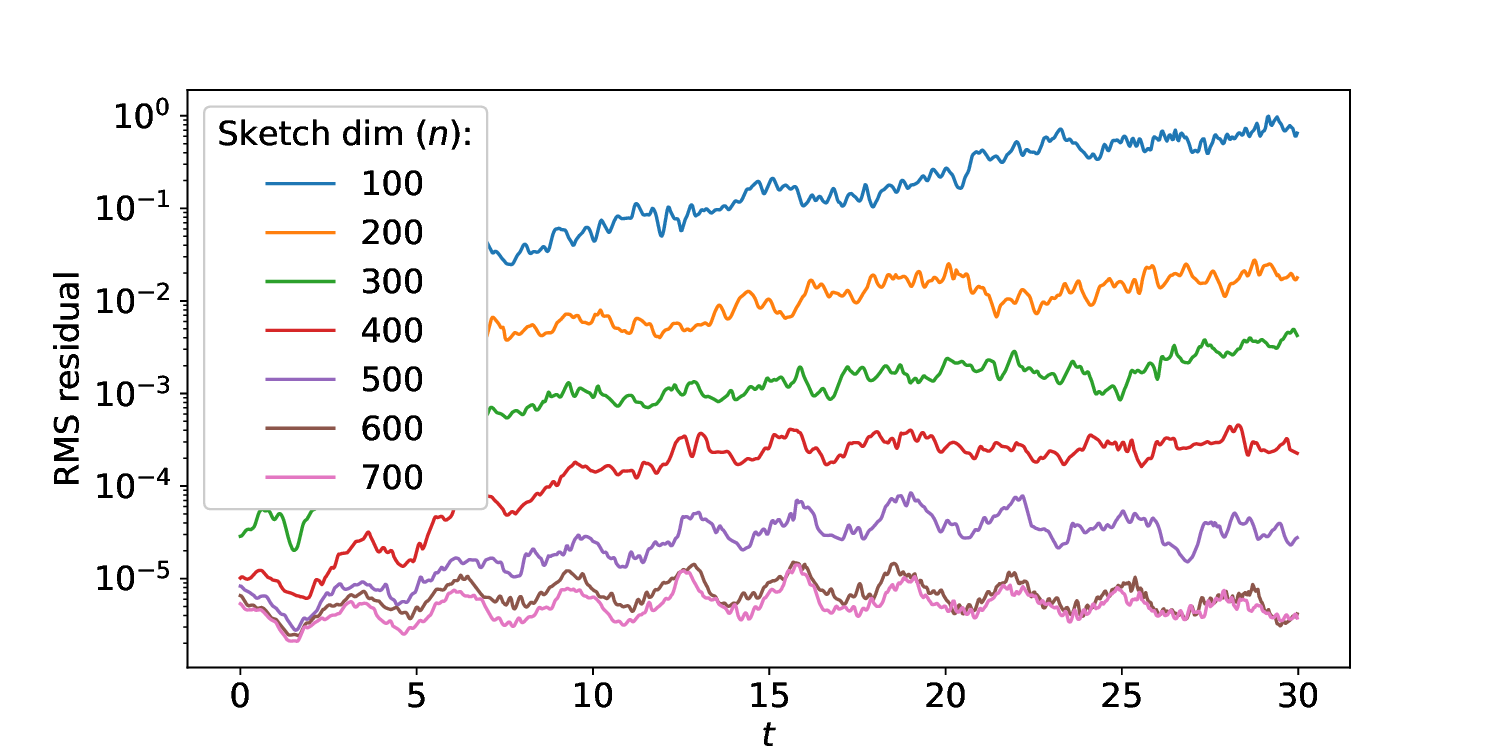}
\par\end{centering}
\caption{RMS residual (\ref{eq:rmseres}) over the MNE trajectory for the example
described in Section \ref{subsec:Langevin-dynamics-with}, plotted
for various choices of the sketch dimension $n$. The regularization
parameter is $\protect\ve=10^{-6}$ for each curve. \label{fig:rank_langevin}}
\end{figure}

\begin{figure}
\begin{centering}
\includegraphics[scale=0.5]{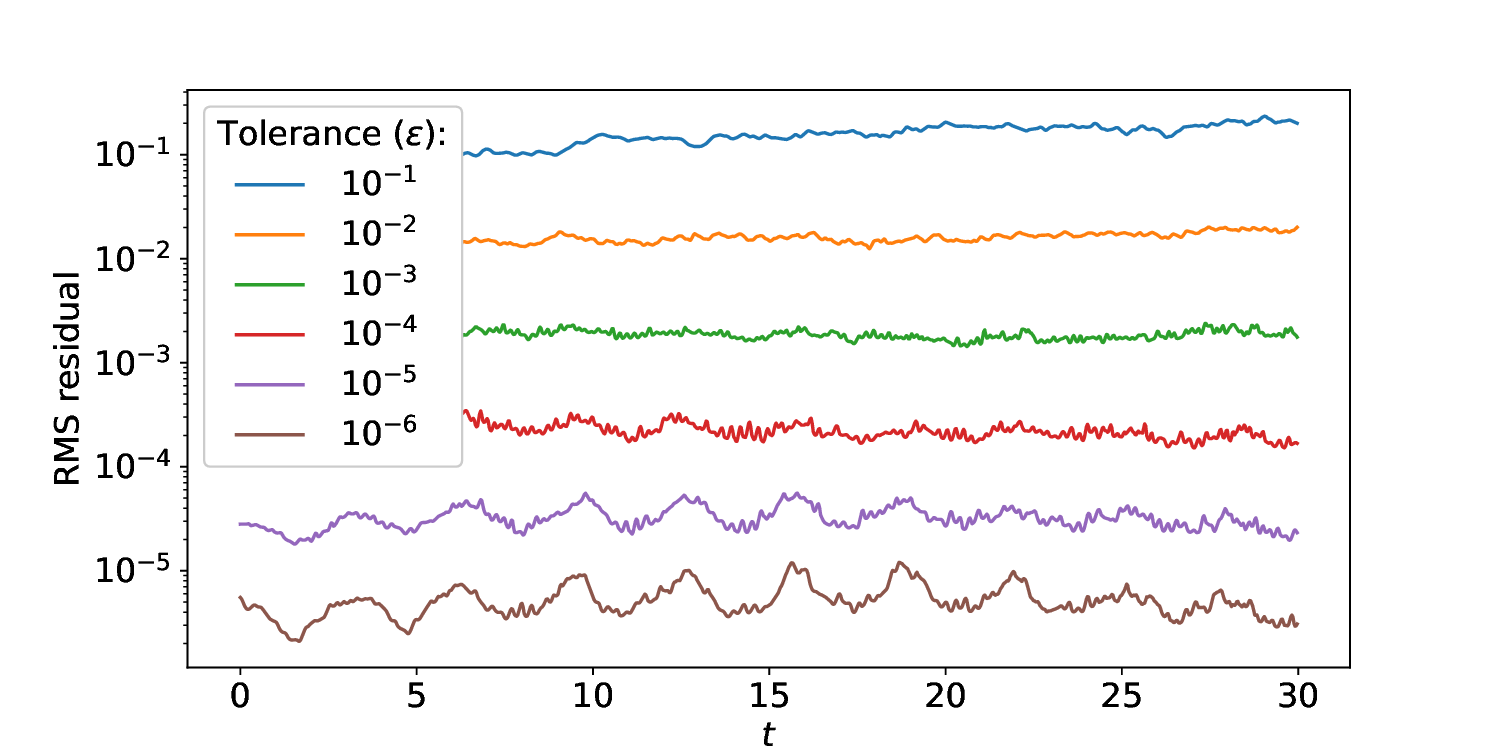}
\par\end{centering}
\caption{RMS residual (\ref{eq:rmseres}) over the MNE trajectory for the example
described in Section \ref{subsec:Langevin-dynamics-with}, plotted
for various choices of the tolerance $\protect\ve$ (cf. (\ref{eq:regparam})).
The sketch dimension is $n=700$ for each curve. \label{fig:tol_langevin}}
\end{figure}

To demonstrate the systematic improvability of the method we plot
the size of the residual for solving (\ref{eq:energyevol}), aggregated
in the root mean square (RMS) sense over the evolving set collocation
points (cf. (\ref{eq:ridge})): 
\begin{equation}
r(t)=\frac{1}{\sqrt{N}}\left\Vert \left[\nabla_{\theta}u_{\theta}(X(t))\right]^{\top}\dot{\theta}(t)-\mathcal{A}_{t}[u_{\theta}](X(t))\right\Vert ,\label{eq:rmseres}
\end{equation}
 for several values of the sketch dimension $n$ and regularization
parameter $\ve$ in Figures \ref{fig:rank_langevin} and \ref{fig:tol_langevin},
respectively. Note that accuracy for our parametrization saturates
by about $n=600$. Moreover, past about $\ve=10^{-6}$, the MNE for
this parametrization fails to generalize outside the collocation points
and the ODE solver becomes unstable. 

Figures \ref{fig:rank_langevin} and \ref{fig:tol_langevin} validate
the systematic improvability of the methodology up to the accuracy
allowed by the parametrization, consistent with prior expectations.
Note that the experiments can be viewed as revealing a numerical rank
for the MNE that is far smaller than the number of parameters.

\subsection{Bayesian inverse problem \label{subsec:Bayesian-inverse-problem}}

We will apply the MNE to sample from the Bayesian posterior distribution
in a fully Bayesian Gaussian process regression (GPR) problem. For
further background, see, e.g., \cite{kielstra2024gradientbaseddeterminantfreeframeworkfully}.
This experiment can be viewed as an application of the techniques
of Section \ref{subsec:Sampling-and-marginalization}, where we choose
$\sigma=\sqrt{2}$ and $\gamma=1$.

We construct a dataset of ordered pairs $(t_{i},y_{i})$, $i=1,\ldots,m$,
as 
\[
y_{i}=f(t_{i})+\epsilon_{i},
\]
 where $\epsilon_{i}\sim\mathcal{N}(0,0.1^{2})$ are i.i.d. (Here
`$t$' does not indicate a time variable in the sense of our diffusion
process, but rather an independent regression variable.) In our experiments,
we draw the $t_{i}$ independently from the uniform distribution over
$[-1,1]$. We also choose $f(t)=\sin(5t)$ and $m=20$.

For our GPR, we consider the exponential kernel 
\[
\Sigma(t_{1},t_{2})=\alpha^{2}e^{-\frac{(t_{1}-t_{2})^{2}}{\rho^{2}}},
\]
 where $\alpha>0$ and $\rho>0$ indicate hyperparameters for the
vertical and horizontal scales of the function. Moreover, we let $\sigma>0$
denote a hyperparameter for the unknown standard deviation of the
noise.

We will write $\alpha=e^{x_{1}}$, $\rho=e^{x_{2}}$, and $\sigma^{x_{3}}$
in terms of hyperparameters $x=(x_{1},x_{2},x_{3})$ to be inferred.
As a prior over $x$ we take the standard normal distribution $p(x)\propto e^{-\Vert x\Vert^{2}/2}$.
Then the posterior distribution for $x$, given our observations $y_{i}$,
is specified by the energy function
\[
u(x)=\frac{1}{2}\log\det\mathbf{K}(x)+\frac{1}{2}y^{\top}\mathbf{K}(x)^{-1}y-\log p(x),
\]
 which we view as specifying our target for sampling.

To fit our surrogate model, we pick the initial collocation points
by running 1000 iterations of the Metropolis-adjsuted Langevin algorithm
(MALA) with step size 0.02 on an ensemble of $N=10^{4}$ walkers initialized
at the origin in $\R^{3}$. For our neural network parametrization
we use the MLP architecture with 2 hidden layers of width 128 and
softplus activations. We fit the surrogate model with Adam \cite{adam}
using standard hyperparameters and learning rate $4\times10^{-4}$.

\begin{figure}
\begin{centering}
\includegraphics[scale=0.5]{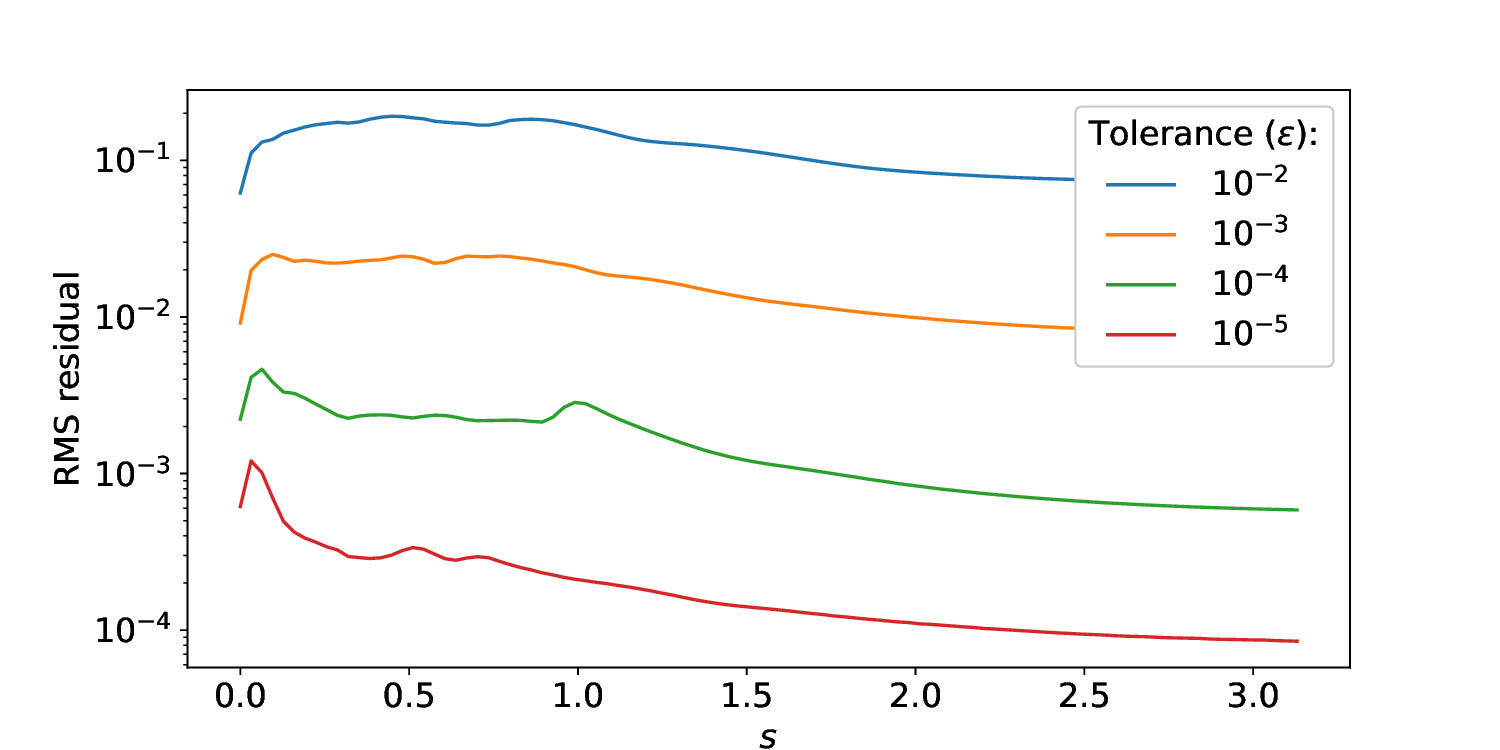}
\par\end{centering}
\caption{RMS residual (\ref{eq:rmseres}) over the MNE trajectory for the example
described in Section \ref{subsec:Bayesian-inverse-problem}, plotted
for various choices of the tolerance $\protect\ve$ (cf. (\ref{eq:regparam})).
The sketch dimension is $n=700$ for each curve. \label{fig:rms_bayesian}}
\end{figure}

We take $s_{\max}=\sqrt{10}$ to be the final time for the MNE. In
Figure \ref{fig:rms_bayesian} we plot the appropriate RMS residual
(cf. (\ref{eq:rmseres})) for this problem as a function of time,
for a fixed sketch dimension $n=700$ and several values of the regularization
parameter $\ve$.

The results once again validate the systematic improvability of the
approach.

\begin{figure}
\begin{centering}
\includegraphics[scale=0.5]{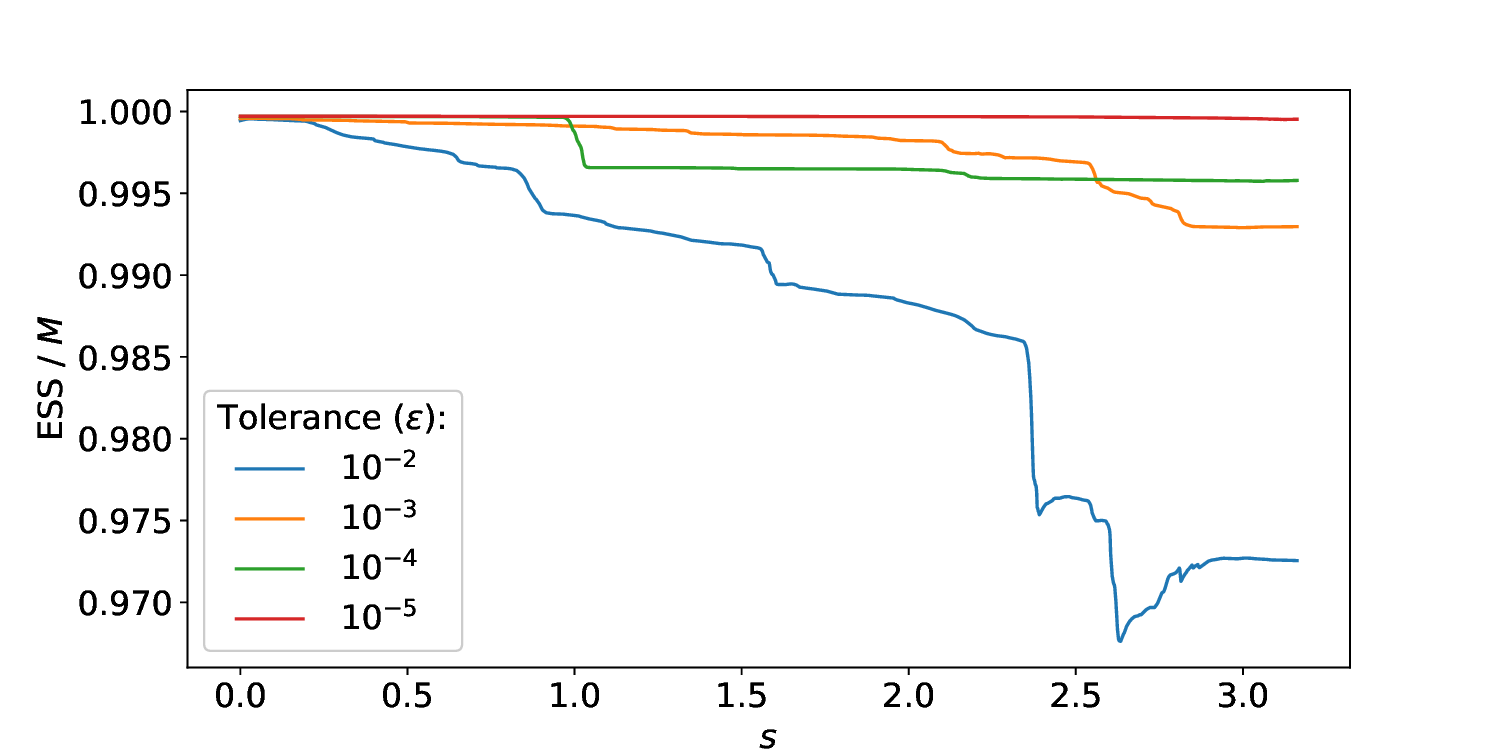}
\par\end{centering}
\caption{ESS per sample for the evolving weighted ensemble of $M=10^{4}$ samples
over the reverse diffusion trajectory for the example described in
Section \ref{subsec:Bayesian-inverse-problem}, plotted for various
choices of the tolerance $\protect\ve$ (cf. (\ref{eq:regparam})).
The sketch dimension is $n=700$ throughout.\label{fig:ess_bayesian}}
\end{figure}

Then we use the learned trajectory $\theta(s)$ from each of these
experiments to draw an unbiased weighted ensemble of $M=10^{4}$ samples
from $e^{-u_{\theta(s_{\max}-s)}}$, evolving with $s\in[0,s_{\max}]$,
following the approach of Section \ref{subsec:Unbiasing}. To solve
the SDE (\ref{eq:revmask-1}), we use the Euler-Maruyama method with
step size $s_{\max}/2000$. The ESS \cite{ESS} per sample is plotted
in Figure \ref{fig:ess_bayesian} over the course of the ensemble
trajectory. In particular, the samples are weighted at $s=0$ according
to the initial step of the unbiasing procedure of Section \ref{subsec:Unbiasing}.
The unbiasing step at the end is ignored here because its impact is
negligible, due to the excellent fit of the surrogate model.

\subsection{Allen-Cahn potential \label{subsec:Allen-Cahn-potential}}

Finally, we consider the high-dimensional distribution induced by
the Allen-Cahn potential for a 1-dimensional field with periodic boundary
condition: 

\[
u(x)=\frac{\beta}{2}\sum_{i=0}^{d-1}\left[\left(\frac{x_{i+1}-x_{i}}{h}\right)^{2}+(x_{i}^{2}-1)^{2}\right].
\]
 The entries of $x$ are zero-indexed with the convention that $x_{d}=x_{0}$.
Here $h>0$ is a spatial grid size parameter and $\beta>0$ is an
inverse temperature. We choose $d=20$, $h=1/20$, and $\beta=0.3$.
This choice of parameters exhibits bimodality in the distribution
$\rho\propto e^{-u}$. In particular the modes are centered at $x=\pm\mathbf{1}$.
This experiment can be viewed as an application of the techniques
of Section \ref{subsec:Sampling-and-marginalization}, where we choose
$\sigma=\sqrt{2}$ and $\gamma=1$.

To fit our surrogate model, we pick the initial collocation points
by running 4000 iterations of MALA with step size $10^{-3}$ on an
ensemble of $N=10^{4}$ walkers initialized according to the standard
normal distribution on $\R^{d}$. We comment that this ensemble is
\emph{far} from fully mixed with respect to the target density $\rho$.
For our neural network parametrization we use the MLP architecture
with 2 hidden layers of width 128 and softplus activations. We fit
the surrogate model with Adam using standard hyperparameters and learning
rate $10^{-4}$.

\begin{figure}
\begin{centering}
\includegraphics[scale=0.5]{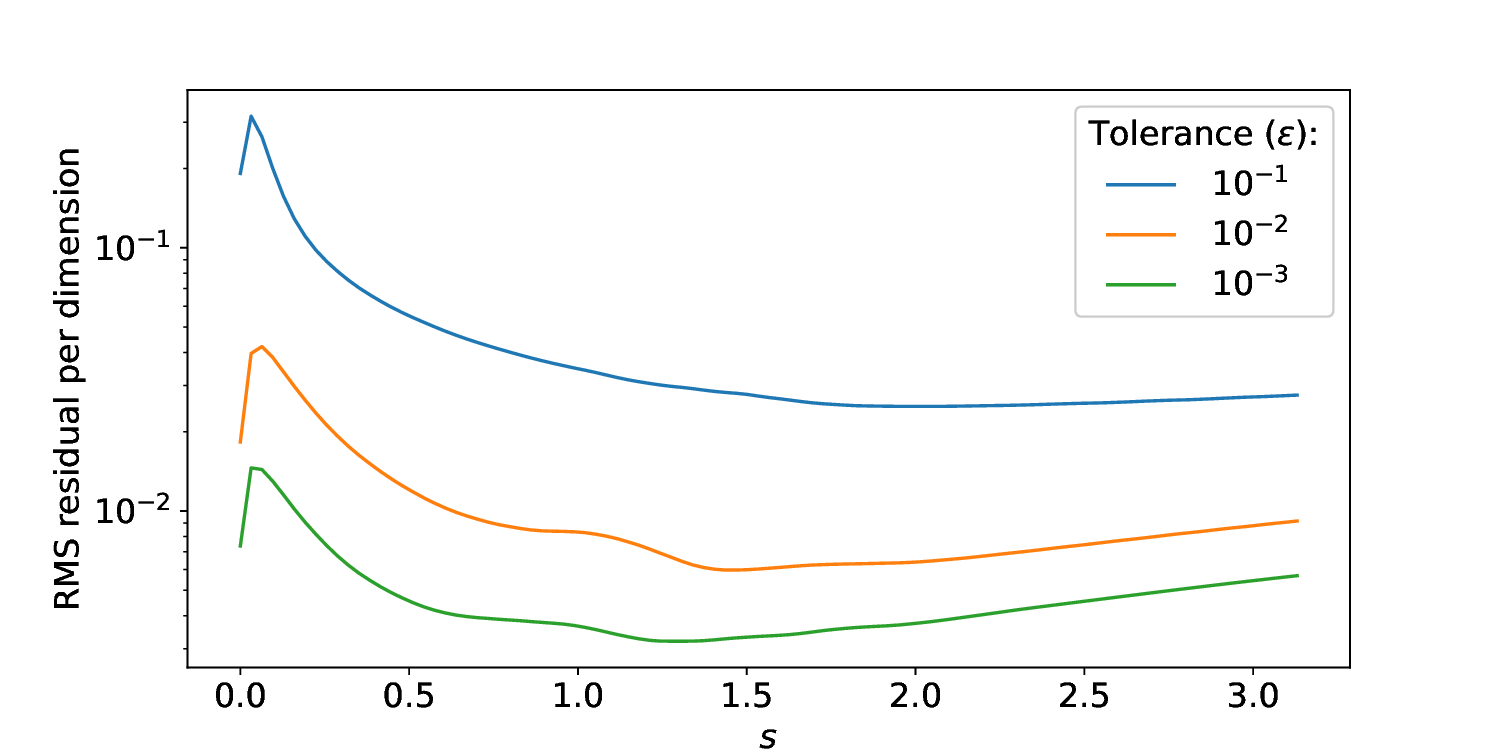}
\par\end{centering}
\caption{RMS residual (\ref{eq:rmseres}) over the MNE trajectory for the example
described in Section \ref{subsec:Allen-Cahn-potential}, plotted for
various choices of the tolerance $\protect\ve$ (cf. (\ref{eq:regparam})).
The sketch dimension is $n=2000$ for each curve. \label{fig:rms_ac}}
\end{figure}

We take $s_{\max}=\sqrt{10}$ to be the final time for the MNE. In
Figure \ref{fig:rms_ac} we plot the appropriate RMS residual (cf.
(\ref{eq:rmseres})) for this problem as a function of time, for a
fixed sketch dimension $n=2000$ and several values of the regularization
parameter $\ve$.

\begin{figure}
\begin{centering}
\includegraphics[scale=0.5]{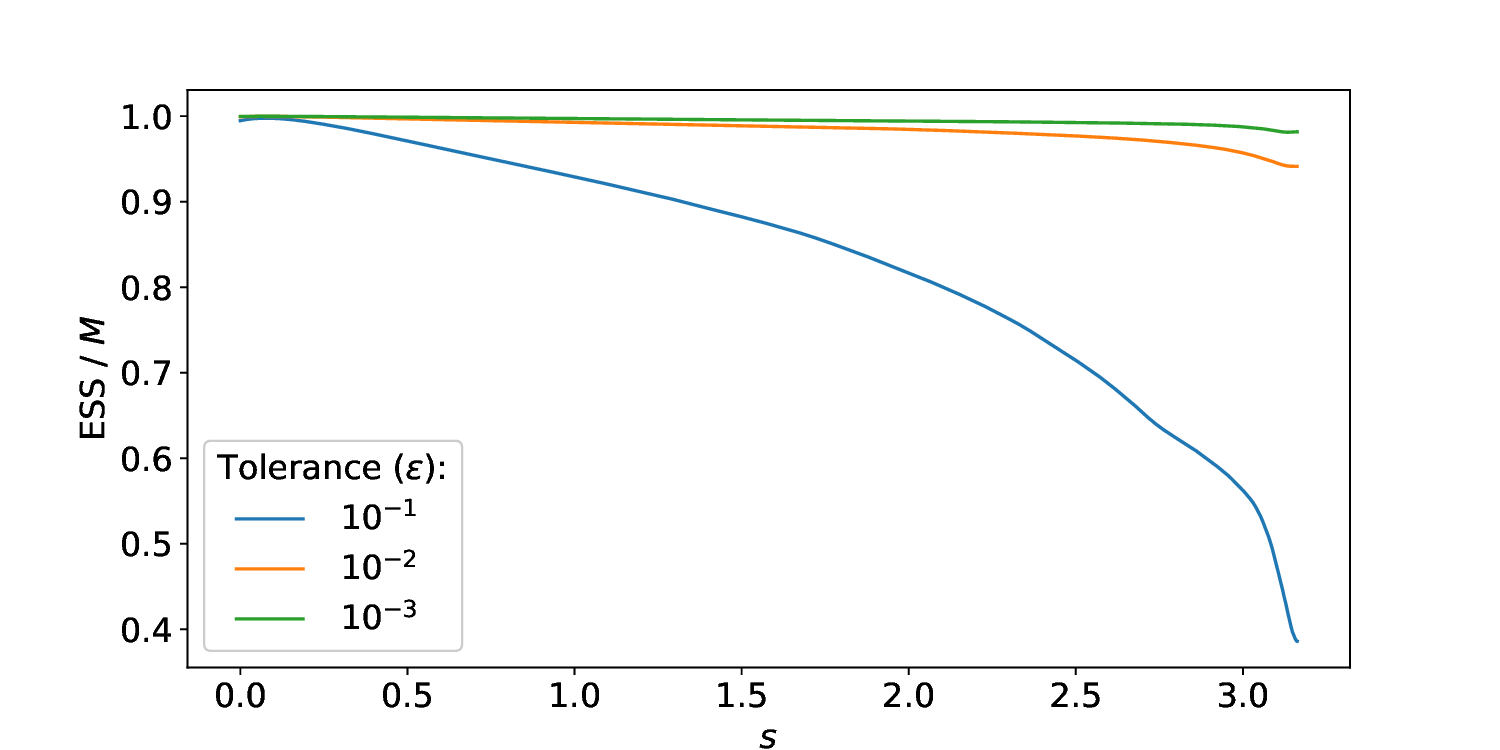}
\par\end{centering}
\caption{ESS per sample for the evolving weighted ensemble of $M=10^{4}$ samples
over the reverse diffusion trajectory for the example described in
Section \ref{subsec:Allen-Cahn-potential}, plotted for various choices
of the tolerance $\protect\ve$ (cf. (\ref{eq:regparam})). The sketch
dimension is $n=2000$ throughout. The final values of the ESS per
sample after the unbiasing step at the end of the reverse diffusion
are $0.386$, $0.937$, and $0.978$ for $\protect\ve=10^{-1}$, $10^{-2}$,
and $10^{-3}$, respectively. \label{fig:ess_ac}}
\end{figure}

Then we use the learned trajectory $\theta(s)$ from each of these
experiments to draw an unbiased weighted ensemble of $M=10^{4}$ samples
from $e^{-u_{\theta(s_{\max}-s)}}$, evolving with $s\in[0,s_{\max}]$,
following the approach of Section \ref{subsec:Unbiasing}. To solve
the SDE (\ref{eq:revmask-1}), we use the Euler-Maruyama method with
step size $s_{\max}/2000$. The ESS per sample is plotted in Figure
\ref{fig:ess_ac} over the course of the ensemble trajectory. In particular,
the samples are weighted at $s=0$ according to the initial step of
the unbiasing procedure of Section \ref{subsec:Unbiasing}. The final
ESS values after the unbiasing step at the end of the reverse diffusion
are reported in the caption of Figure \ref{fig:ess_ac}.

\begin{figure}
\begin{centering}
\includegraphics[scale=0.5]{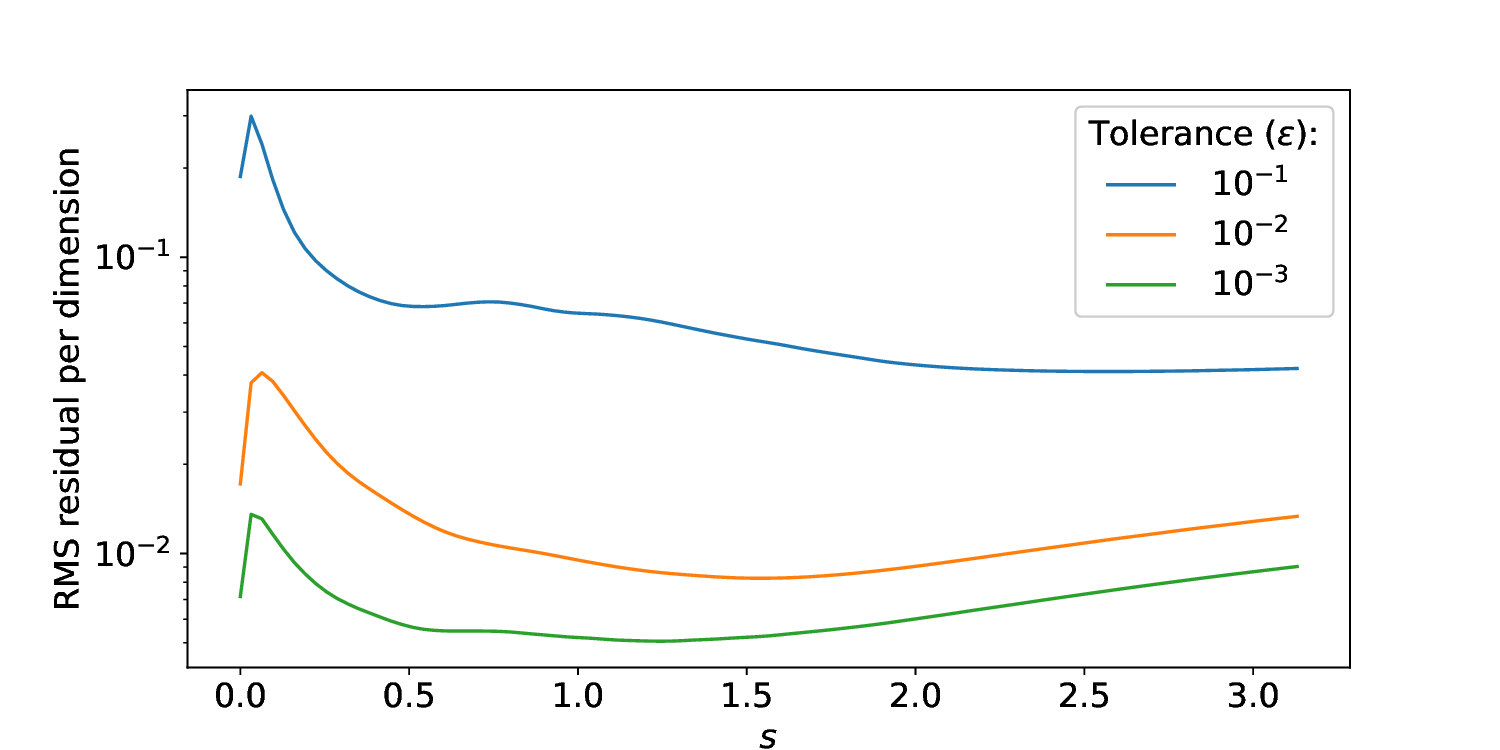}
\par\end{centering}
\caption{RMS residual (\ref{eq:rmseres}) over the MNE trajectory for the example
described in Section \ref{subsec:Allen-Cahn-potential} with $\mathcal{S}=\{1,\ldots,d-1\}$
(i.e., leaving out the zeroth index), plotted for various choices
of the tolerance $\protect\ve$ (cf. (\ref{eq:regparam})). The sketch
dimension is $n=2000$ for each curve. The final values of the ESS
per sample after the unbiasing step at the end of the reverse diffusion
are $0.149$, $0.923$, and $0.967$ for $\protect\ve=10^{-1}$, $10^{-2}$,
and $10^{-3}$, respectively. \label{fig:rms_ac_cond}}
\end{figure}

Finally, we apply the MNE to conditional sampling, choosing $\mathcal{S}=\{1,\ldots,d-1\}$
(i.e., leaving out the zeroth index) in the notation of Section \ref{subsec:Sampling-and-marginalization}
and conditioning on $x_{0}=1$, which essentially confines us to a
single mode of the bimodal distribution. The setup for this experiment
is the same as in the preceding discussion. In Figure \ref{fig:rms_ac_cond}
we plot the appropriate RMS residual (cf. (\ref{eq:rmseres})) for
this problem as a function of time.

The ESS plot analogous to Figure \ref{fig:ess_ac} is qualitatively
quite similar for this case of conditional sampling. Thus we omit
the plot, but we include the final ESS values in the caption of Figure
\ref{fig:rms_ac_cond}, which are quite similar to those obtained
in the full sampling case.

\bibliographystyle{plain}
\bibliography{mne}

\end{document}